\documentclass[11pt]{article}
\usepackage[a4paper,margin=1in]{geometry}
\usepackage{amsmath,amssymb,amsfonts,amsthm,mathtools}
\usepackage{enumitem}
\usepackage{comment}
\usepackage{needspace}
\usepackage{microtype}
\usepackage{lmodern}
\usepackage{color}
\usepackage[T1]{fontenc}
\setlist[enumerate]{itemsep=2pt,topsep=4pt}
\allowdisplaybreaks[1]
\usepackage[hidelinks]{hyperref}
\usepackage{aliascnt}
\usepackage[nameinlink,capitalize]{cleveref}
\newtheorem{thm}{Theorem}[section]
\newaliascnt{lemma}{thm}
\newtheorem{lem}[lemma]{Lemma}
\aliascntresetthe{lemma}
\newaliascnt{proposition}{thm}
\newtheorem{prop}[proposition]{Proposition}
\aliascntresetthe{proposition}
\newaliascnt{corollary}{thm}
\newtheorem{corollary}[corollary]{Corollary}
\aliascntresetthe{corollary}
\newaliascnt{claim}{thm}
\newtheorem{clm}[claim]{Claim}
\aliascntresetthe{claim}
\theoremstyle{definition}
\newaliascnt{definition}{thm}
\newtheorem{defi}[definition]{Definition}
\aliascntresetthe{definition}
\newaliascnt{remark}{thm}

\aliascntresetthe{remark}

\allowdisplaybreaks 
\title{P\'osa-type resilience for Hamilton Berge cycles\\ in random hypergraphs\thanks{The work is supported by NSFC(12171088, 11871015) and CGLSTDFFJ(2023L3003). E-mail: 1519036064@qq.com (L. Zhang),  elian1425@fzu.edu.cn (A. Chen), qliu@fzu.edu.cn (Q. Liu). Mathematics Subject Classification [2020]:  05C45, 05C70, 05C80. } }
\author{Liping Zhang$^{a}$,
	Ailian Chen$^{b}$\thanks{Corresponding author.}, 
	Qinghai Liu$^{a,c}$\\[1mm]
{\small\it $^a$Center for Discrete Mathematics, Fuzhou University,
	Fuzhou, Fujian, 350108 China}\\
{\small\it $^b$School of Mathematics and Statistics, Fuzhou University, Fuzhou, Fujian, 350108 China}\\
{\small\it $^c$Key Laboratory for Operations Research and Cybernetics of Fujian Universities}}
\date{}
 
\begin{document}
\maketitle
		\begin{abstract}
		We prove a degree-sequence resilience theorem for Hamilton Berge cycles
		in random uniform hypergraphs. For every fixed integer $k\ge3$ and
		every $\gamma>0$, there exists $C=C(k,\gamma)>0$ such that if
		$C\log n/n^{k-1}\le p\le1$, then asymptotically almost surely every
		spanning subhypergraph of $H^{(k)}(n,p)$ with increasing vertex-degree
		sequence $g_1\le\cdots\le g_n$ satisfying
		\[
		g_i\ge\left(\left(\frac{i}{n}\right)^{k-1}+\gamma\right)
		\binom{n-1}{k-1}p
		\qquad(1\le i<n/2)
		\]
		contains a Hamilton Berge cycle. 
		The term $(i/n)^{k-1}$ agrees with
		the normalized asymptotic degree bound in Salia's deterministic
		P\'osa-type theorem.
		Our result gives a hypergraph analogue of the random-graph theorem
		of Condon et al. and strengthens the sharp Dirac-type resilience
		theorem of Chen and Zhang at the same logarithmic probability order.

		The deletion-resilience formulation has optimal probability order up
		to constant factors. The coefficient $1$ of $(i/n)^{k-1}$ cannot be
		replaced by a smaller constant. The proof controls adaptive rotations
		and representative-edge incidences uniformly. It also gives an endpoint-growth estimate under nonuniform
		deletion bounds.

	\end{abstract}
	
	\section{Introduction} A \emph{Hamilton cycle} of a graph is a cycle passing through every vertex exactly once, and a graph is called \emph{Hamiltonian} if it contains a Hamilton cycle. Hamiltonicity is one of the classical topics in graph theory. Since deciding whether a given graph contains a Hamilton cycle is NP-complete \cite{Karp}, it is natural and important to seek sufficient conditions which guarantee Hamiltonicity.
	Degree conditions for Hamiltonicity connect extremal graph theory with
	probabilistic combinatorics. Dirac's theorem~\cite{Dirac} states that
	every graph on $n\ge3$ vertices with minimum degree at least $n/2$ is
	Hamiltonian. P\'osa~\cite{Posa1962} strengthened this by showing that an
	increasing degree sequence $d_1\le\cdots\le d_n$ guarantees a Hamilton
	cycle whenever $d_i\ge i+1$ for every integer $1\le i<n/2$.
	Chv\'atal~\cite{Chvatal} subsequently obtained the more general
	condition that, for every $i<n/2$, either $d_i\ge i+1$ or
	$d_{n-i}\ge n-i$.
	
	Resilience measures how robust such conclusions remain under arbitrary edge deletions. Its systematic study was initiated by Sudakov and
	Vu~\cite{SudakovVu}. The binomial random graph $G(n,p)$ is obtained from the complete graph on $n$ vertices by letting each edge appear independently with probability $p$. We say that an event
	holds \emph{asymptotically almost surely}, or a.a.s., if its probability
	tends to one as $n\to\infty$. Lee and Sudakov~\cite{LeeSudakov} proved
	that, for every fixed $\gamma>0$ and sufficiently large $C$, if
	$p\ge C\log n/n$, then a.a.s. every spanning subgraph of $G(n,p)$
	with minimum degree at least $(1/2+\gamma)np$ is Hamiltonian.
	Condon et al.~\cite{CondonEtAl} obtained the corresponding
	P\'osa-resilience theorem at the same probability order.
	Fix $\gamma>0$ and let $C$ be sufficiently large.
	For $p\ge C\log n/n$, a.a.s. every spanning subgraph is Hamiltonian
	if its increasing degree sequence satisfies
	$d_i\ge(i/n+\gamma)np$ for $1\le i<n/2$.
	Their theorem uses the decreasing degrees of the deleted subgraph.
	We use this formulation below.
	
	A \emph{$k$-uniform hypergraph}, or \emph{$k$-graph}, is a pair
	$H=(V,E)$ with $E\subseteq\binom Vk$. The random $k$-graph
	$H^{(k)}(n,p)$ has vertex set $[n]=\{1,\ldots,n\}$ and contains each
	possible hyperedge independently with probability $p$.
	A \emph{Berge cycle} of length $t\ge2$ consists of distinct \emph{core
		vertices} $v_1,\ldots,v_t$ and distinct hyperedges $e_1,\ldots,e_t$
	such that $\{v_i,v_{i+1}\}\subseteq e_i$ for all $ i\in[t], $ where $v_{t+1}=v_1$~\cite{Berge}. It is a \emph{Hamilton Berge cycle} if its core vertex
	set is the entire vertex set of the hypergraph. We call a hypergraph
	containing such a cycle \emph{Berge Hamiltonian}.  The study of Berge cycles in hypergraphs has a long history.
	
	There is a substantial deterministic theory of degree conditions for
	Berge cycles. Early Dirac-type results go back to Bermond, Germa,
	Heydemann and Sotteau~\cite{BermondEtAl}. F\"uredi, Kostochka and
	Luo~\cite{FurediKostochkaLuo} obtained sharp minimum-degree conditions
	for long Berge cycles in nonuniform hypergraphs, while Kostochka, Luo
	and McCourt~\cite{KostochkaLuoMcCourt} proved exact Dirac-type bounds
	in uniform hypergraphs. Salia~\cite[Theorem~5]{Salia} proved an
	exact P\'osa-type degree-sequence theorem for $k\ge3$ and $n>2k$.
	It requires $d_i>i$ for $1\le i<k$.
	It also requires  $d_i>\binom{i}{k-1}$ for $k\le i<n/2$.
	For even $n$, it also requires
	$d_{n/2-1}>\binom{n/2-1}{k-1}+1$.
	These conditions guarantee a Hamilton Berge cycle.

	For random $k$-graphs, Clemens, Ehrenm\"uller and
	Person~\cite{ClemensEhrenmullerPerson} proved that for every fixed
	$k\ge3$ and $\gamma>0$, if
	$p\ge(\log n)^{17k}/n^{k-1}$, then a.a.s. every spanning subhypergraph
	with minimum vertex degree at least
	$(2^{-(k-1)}+\gamma)\binom{n-1}{k-1}p$ is Berge Hamiltonian.
	Chen and Zhang~\cite{ChenZhang} obtained this sharp minimum-degree
	constant for $k\ge4$ and 
	$C\log n/n^{k-1}\le p\le1$, with $C=C(k,\gamma)$. Their argument
	also recovers the $k=3$ case. An earlier $3$-uniform result appears
	in~\cite{ChenZhang3}. Thus a natural next question is whether a
	nonuniform degree-sequence condition can replace the minimum-degree
	condition at the same logarithmic probability scale.

	The threshold for the random hypergraph itself to be Berge Hamiltonian is
	known more precisely. Bal, Berkowitz, Devlin and
	Schacht~\cite{BalBerkowitzDevlinSchacht} proved that, in the random
	hypergraph process, a Hamilton Berge cycle a.a.s. appears as soon as
	the minimum vertex degree becomes two. In particular, if
	\[
	p=\frac{(k-1)!(\log n+\log\log n+c_n)}{n^{k-1}},
	\qquad c_n\to c\in\mathbb R,
	\]
	the probability of Berge Hamiltonicity tends to $\exp(-e^{-c})$. This is a statement about the
	undeleted random
	$k$-graphs, rather than  resilience. Related resilience
	questions for other notions of hypergraph Hamilton cycles have also
	been studied, including tight and loose cycles;
	see~\cite{AllenParczykPfenninger,PetrovaTrujic,KeevashEtAl,NenadovSkoric}.

	We prove a P\'osa-type resilience theorem for Hamilton Berge cycles
	for every fixed $k\ge3$.
	Before stating the theorem, we introduce the deletion formulation.
	Unless specified otherwise, for $H^{(k)}(n,p)$ we use 
	$D:=\binom n{k-1}p$.
	
	\begin{defi}[$\varepsilon$-P\'osa-resilience]\label{def-resi}
		Fix a vertex set $V=[n]$, an integer $k\ge3$, and parameters
		$D>0$ and $\varepsilon>0$. Let  
		$
		\mathcal H_{D}^{\varepsilon}$
		denote the family of all $k$-graphs $H'$ on this fixed vertex set $V$
		for which there exists an ordering
		$v_1,\ldots,v_n$
		of the vertices satisfying
		$d_{H'}(v_1)\ge d_{H'}(v_2)\ge \cdots \ge d_{H'}(v_n)$
		and
		\begin{align}\label{1.1}
		d_{H'}(v_i)
		\le \left(1-\left(\frac{i}{n}\right)^{k-1}
		-4k\varepsilon\right)D
		\qquad\text{for every integer }1\le i<\frac n2.
		\end{align}
		For a $k$-graph $H$ on the same vertex set $V$, define $
		\mathcal H_{D}^{\varepsilon}(H)
		=\left\{
		H'\in \mathcal H_{D}^{\varepsilon}
		:H'\subseteq H
		\right\}.
		$
		Here $H-H'$ has vertex set $V$ and edge set $E(H)\setminus E(H')$
		(all deletions are edge deletions). We say that $H$ is
		$\varepsilon$-P\'osa-resilient with respect to a property $\mathcal P$ if
		$H-H'
		\in \mathcal P$
		for every
		$H'\in \mathcal H_{D}^{\varepsilon}(H)$.
	\end{defi}
	In particular, every deleted degree is at most
	$(1-n^{-(k-1)}-4k\varepsilon)D$.
	Our main result is the following.
	
	\begin{thm}\label{thm-min}
		For every fixed integer \(k\ge 3\) and every \(\varepsilon>0\), there exists a constant \(C>0\) such that if
		$C\frac{\log n}{n^{k-1}}\le p\le1,$
		then a.a.s. the random $k$-graph
		$
		H^{(k)}(n,p)$
		is \(\varepsilon\)-P\'osa-resilient with respect to Berge Hamiltonicity.
	\end{thm}
	
	The probability order in the deletion-resilience
	Theorem~\ref{thm-min} is also optimal up to constant factors.
	For sufficiently small fixed $\varepsilon>0$, the empty deletion is
	admissible.
	The host itself must therefore be Berge Hamiltonian.
	If $0<p\ll\log n/n^{k-1}$, it a.a.s. has an isolated vertex;
	see~\cite{BalBerkowitzDevlinSchacht}.
	This rules out a smaller probability order for deletion-resilience.
	We do not determine the optimal constant in that probability bound. The following two results indicate that the P\'osa-type degree-sequence condition is asymptotically best possible.
	
	The deletion formulation yields a condition on the increasing degree
	sequence of the retained hypergraph.
	Write $D_0:=\binom{n-1}{k-1}p$.
	The $D$ in Definition~\ref{def-resi} satisfies
	\begin{equation}\label{eq:normalization-D-D0}
	D=\binom n{k-1}p=\frac{n}{n-k+1}D_0\ge D_0.
	\end{equation}
	
	\begin{corollary}[P\'osa-type consequence]\label{cor:retained-sequence}
		For every fixed integer $k\ge3$ and every $\gamma>0$, there exists
		$C=C(k,\gamma)>0$ such that the following holds whenever
		$C\log n/n^{k-1}\le p\le1$.
		Asymptotically almost surely, every spanning subhypergraph
		$G\subseteq H^{(k)}(n,p)$ whose increasing vertex-degree sequence
		$g_1\le\cdots\le g_n$ satisfies
		\begin{equation}\label{eq:retained-degree-condition}
		g_i\ge\left(\left(\frac{i}{n}\right)^{k-1}+\gamma\right)D_0
		\qquad\text{for every integer }1\le i<n/2
		\end{equation}
		is Berge Hamiltonian.
	\end{corollary}
	
	\begin{proof}
		It suffices to prove the assertion for $0<\gamma\le1/4$, since
		the hypothesis becomes stronger as $\gamma$ increases.
		Take $\varepsilon=\gamma/(8k)$ and choose $C$ sufficiently large.
		Let $H=H^{(k)}(n,p)$.
		Each  degree of $H$ is binomial with mean $D_0=\Omega_k(C\log n)$.
		Chernoff's inequality (\cref{thm-cher}) and a union bound give
		$\Delta(H)\le(1+\gamma/4)D_0$ a.a.s.
		Condition on this event and on the event in Theorem~\ref{thm-min}
		with parameter $\varepsilon$.
		
		Let $G\subseteq H$ satisfy~\eqref{eq:retained-degree-condition}, and
		let $F:=H-G$. Write $f_1\ge\cdots\ge f_n$ for the decreasing
		degree sequence of $F$. Fix an integer $1\le i<n/2$. At least
		$n-i+1$ vertices $v$ satisfy $d_G(v)\ge g_i$.
		For each such vertex,
		$d_F(v)\le(1-(i/n)^{k-1}-3\gamma/4)D_0$.
		At most $i-1$ vertices can exceed this bound. Hence
		\[
		\begin{aligned}
		f_i
		&\le\left(1-\left(\frac{i}{n}\right)^{k-1}-\frac{3\gamma}{4}\right)D_0\le\left(1-\left(\frac{i}{n}\right)^{k-1}-4k\varepsilon\right)D_0\\
		&\le\left(1-\left(\frac{i}{n}\right)^{k-1}-4k\varepsilon\right)D.
		\end{aligned}
		\]
		The last inequality uses~\eqref{eq:normalization-D-D0} and the
		positivity of its coefficient. Since $i<n/2$, $k\ge3$, and
		$\gamma\le1/4$, that coefficient is at least
		$1-1/4-1/8>0$. Thus $F\in\mathcal H_D^\varepsilon(H)$.
		Theorem~\ref{thm-min} gives that $G=H-F$ is Berge Hamiltonian.
		These events depend only on  $H$ and hold 
		for every spanning subhypergraph $G$ under consideration.
	\end{proof}
	
	Our degree condition has the same leading term as the condition in~\cite{Salia}.
	Indeed, $\binom{i}{k-1}/\binom{n-1}{k-1}
	=(i/n)^{k-1}+O_k(n^{-1})$  for $1\le i<n/2$.
	Here the binomial coefficient is zero when $i<k-1$.
	We preserve this condition after scaling by $p$.

	The two theorems address different settings.
	At $p=1$, our condition still contains the fixed term 
	$\gamma\binom{n-1}{k-1}$.
	It therefore does not recover Salia's exact bounds.
	When $p$ tends to zero, our degree thresholds at linear indices are
	smaller by a factor of order $p$.
	Salia's unscaled degree bounds do not imply this sparse resilience result.

	At the middle index the power condition reaches the Dirac constant
	$2^{-(k-1)}$, while at earlier indices it is strictly below the linear
	condition $(2^{-(k-2)}i/n+\gamma)\binom{n-1}{k-1}p$. This is a hypergraph analogue of the
	degree-sequence strengthening obtained by Condon et al.~\cite{CondonEtAl}.

	Substitute $k=2$
	and $D=np$ in~\eqref{1.1}.
	This gives $d_{H'}(v_i)\le(n-i)p-8\varepsilon np$ for $i<n/2$.
	After reparametrizing the positive constant, this agrees with the
	P\'osa-resilience condition of Condon et al.~\cite{CondonEtAl}.
	\begin{prop}
		\label{prop:uniform-coefficient-sharp}
		Fix $k\ge3$ and $0\le c<1$.
		There are constants $\gamma>0$ and $C_0>0$ with the following property.
		If $C_0\log n/n^{k-1}\le p\le1$, then a.a.s.
		$H^{(k)}(n,p)$ contains a spanning subhypergraph $G$ with no Hamilton
		Berge cycle.
		Its increasing degree sequence satisfies
		$g_i\ge\bigl(c(i/n)^{k-1}+\gamma\bigr)D_0$ for every integer
		$1\le i<n/2$.
	\end{prop}
	
	\begin{proof}
		Take $\gamma=(1-c)2^{-(k+1)}$.
		Apply \cite[Proposition~1.2]{ChenZhang} with error $\gamma$.
		For sufficiently large $C_0$, a.a.s. it supplies a spanning $G$ with
		no Hamilton Berge cycle and $\delta(G)\ge(2^{-(k-1)}-\gamma)D_0$.
		For $i<n/2$, we have
		$c(i/n)^{k-1}+\gamma\le c\,2^{-(k-1)}+\gamma
		=2^{-(k-1)}-3\gamma<2^{-(k-1)}-\gamma$.
		Thus every required degree inequality holds.
		The margin $\gamma$ is positive and independent of $n,p$.
	\end{proof}
	
	This proves tightness of the common leading coefficient.
	Sharpness at each individual index requires separate obstructions.
	Salia's deterministic theorem has such examples~\cite[Section~1.3]{Salia}.
	The same balanced cut also proves asymptotic optimality of the Dirac
	constant $2^{-(k-1)}$; see \cite[Proposition~1.2]{ChenZhang}.

	\begin{corollary}[Dirac-type consequence]\label{cor:dirac-consequence}
		For every fixed integer $k\ge3$ and every $\gamma>0$, there exists
		$C=C(k,\gamma)>0$ such that if
		$C\log n/n^{k-1}\le p\le1$, then a.a.s. every spanning subhypergraph $G\subseteq H^{(k)}(n,p)$ satisfying
		\[
		\delta(G)\ge\left(2^{-(k-1)}+\gamma\right)
		\binom{n-1}{k-1}p
		\]
		is Berge Hamiltonian.
	\end{corollary}
	
	\begin{proof}
		For $i<n/2$, we have $(i/n)^{k-1}\le2^{-(k-1)}$.
		Thus every $g_i\ge\delta(G)$ satisfies the condition in
		Corollary~\ref{cor:retained-sequence}.
	\end{proof}
	
	Corollary~\ref{cor:dirac-consequence} recovers the minimum-degree
	conclusion of Chen and Zhang~\cite{ChenZhang}, with a fixed positive
	constant $\gamma$ and at the same probability order.
	For $0<x<1/2$ and $k\ge3$, we have $x^{k-1}<2^{-(k-2)}x$.
	Thus Corollary~\ref{cor:retained-sequence} also yields the sufficient
	linear condition $g_i\ge(2^{-(k-2)}i/n+\gamma)D_0$. 
	
	\noindent\textbf{Proof strategy and the role of the degree sequence.}
	The proof follows the rotation--extension and sparse-complementing
	scheme used for random graphs by Lee and Sudakov~\cite{LeeSudakov}.
	The degree-sensitive endpoint estimate~\cite[Lemma~18]{CondonEtAl}
	allows the deletion-degree indices to drive endpoint growth. In the
	hypergraph setting, the corresponding growth estimate involves the
	$(k-1)$st power of the endpoint density. Both unused-edge and
	used-edge rotations are needed, the latter already appears in
	Salia's work~\cite[Section~2]{Salia}.
	
	Our argument also builds on the uniform rotation estimates of
	Chen and Zhang~\cite{ChenZhang}. In particular, the weighted incidence
	counting and the decomposition into deleted edges and path
	representatives use the same underlying mechanism. For an endpoint
	family whose representatives have rotation depth at most $t$, the
	representative incidences inherited from the common path contribute
	at most $kn$, and the additional ones at most $t|X|$.
	The present proof makes the estimate uniform over endpoint sets of
	any fixed positive linear size and over the target endpoint density,
	rather than only the terminal half-density threshold needed for a
	uniform minimum-degree bound. The block encodings in
	Lemma~\ref{lem:degree-endpoint-growth} also allow separately chosen
	short rotation sequences without requiring a common retained rotation tree.

	Section~2 develops the probabilistic estimates and the Berge-path
	rotation framework. In Section~3.1 we prove the endpoint-growth lemma
	and obtain a sparse spanning subhypergraph with property
	$RE(1/2+\varepsilon)$ after every deletion from $\mathcal H_D^{2\varepsilon}$.
	Section~3.2 proves that the remaining hypergraph complements each such
	sparse subhypergraph.
	The closure result in \cite[Proposition~2.5]{ChenZhang} then gives a
	Hamilton Berge cycle. The appendix contains the sampling argument.
	
	\section{Preliminaries  }
	We first give notation and probabilistic estimates.
	We then develop the Berge-path rotations used in the proof of
	Theorem~\ref{thm-min}.

	\subsection{Basic notation and probabilistic tools}
	
	\noindent\textbf{Notation.}
	For a $k$-graph $H=(V,E)$, write $V(H)=V$, $E(H)=E$ and $e(H)=|E|$.
	For $U\subseteq V$, let $H[U]$ be the subhypergraph induced by $U$;
	its edges are precisely the edges of $H$ contained in $U$.
	The degree $d_H(v)$ is the number of edges containing $v$.
	Write $\delta(H)$ and $\Delta(H)$ for the minimum and maximum vertex degrees.
	The neighbourhood of $v$ is
	$N_H(v)=\{u\in V\setminus\{v\}:\{u,v\}\subseteq e
	\text{ for some }e\in E(H)\}$.
	Thus $d_H(v)$ counts hyperedges, whereas $|N_H(v)|$ counts distinct vertices.
	Different edges through $v$ may contain the same neighbour.
	In particular, $|N_H(v)|\le(k-1)d_H(v)$, with equality only when those
	edges are pairwise disjoint outside $v$.
	
	All logarithms are natural, and all limits are as $n\to\infty$.
	The edge probability $p=p(n)$ may vary with $n$.
	The integer $k$ and parameters declared fixed remain independent of $n$.
	We use the standard $O$, $o$, $\Omega$ and $\Theta$ notation.
	In particular, $o(1)$ denotes a quantity tending to zero.
	Subscripts record the fixed parameters on which implicit constants
	or rates of convergence may depend.

	We use the following standard Chernoff bounds; see~\cite{AlonSpencer}.
	For the weighted lower-tail bound, see also~\cite[Theorem~2.6]{ChenZhang}.
	The weighted bounds apply when one hyperedge contributes to several
	incidences in the same count.
	
	\begin{thm}[Chernoff's inequality]\label{thm-cher}
		Let $X\sim\mathrm{Bi}(N,p)$ and let $\mu=\mathbb E(X)=Np>0$.
		For every $\varepsilon>0$ and every $\lambda\ge\mu$, respectively,
		\[
		\Pr(|X-\mu|>\varepsilon\mu)
		\le2\exp\!\left(-\frac{\varepsilon^2\mu}{2+\varepsilon}\right),\quad
		\Pr(X\ge\lambda)
		\le\left(\frac{e\mu}{\lambda}\right)^\lambda.
		\]
		More generally, let $X=\sum_{j=1}^m a_j\xi_j$, where the $\xi_j$ are
		independent Bernoulli random variables, possibly with different means.
		Suppose that the weights are deterministic and satisfy $0\le a_j\le b$
		for some $b>0$. Let $\mu=\mathbb E(X)>0$. Then we have
		\[
		\begin{aligned}
		\Pr(X\le\mu-t)
		&\le\exp\!\left(-\frac{t^2}{2b\mu}\right)
		&&(0\le t\le\mu),\\
		\Pr(X\ge\mu+t)
		&\le\exp\!\left(-\frac{t^2}{2b(\mu+t/3)}\right)
		&&(t\ge0).
		\end{aligned}
		\]
	\end{thm}
	For sets $A,B\subseteq[n]$,
	define
	\[
	J_H(A,B):=\sum_{x\in A}\sum_{y\in B\setminus\{x\}}d_H(x,y),
	\]
	where $d_H(x,y)$ counts edges containing both $x$ and $y$.
	
	\begin{lem}[Adapted from {\cite[Lemmas~2.7 and~2.8]{ChenZhang}}]
		\label{lem-adjac}
		Fix $k\ge3$, $\varepsilon>0$ and $\rho>0$.
		There are constants $C=C(k,\varepsilon,\rho)>0$ and
		$\lambda=\lambda(k,\rho)>0$ such that, if
		$C\log n/n^{k-1}\le p\le(\log n)^2/n^{k-1}$, then a.a.s.
		$H=H^{(k)}(n,p)$ satisfies the following properties.
		\begin{enumerate}[label=\textup{(\roman*)}]
			\item Every vertex $v$ satisfies
			$(1-\varepsilon)D\le d_H(v)\le(1+\varepsilon)D$.
			
			\item For all $X,W\subseteq[n]$ with
			$1\le |X|\le\lfloor n/((\log n)^{1/4}D_0)\rfloor$
			and $|W|\le2\lambda |X|D_0$, we have
			$J_H(X,W)<{\rho\over 2}|X|D_0$.
		\end{enumerate}
	\end{lem}
	
	\begin{proof}
		For \textup{(i)}, apply~\cite[Lemma~2.7(i)]{ChenZhang}
		with a smaller $\varepsilon$.
		Then use $D/D_0=1+O_k(n^{-1})$.
		
		For \textup{(ii)}, the count in~\cite[Lemma~2.8]{ChenZhang}
		counts an incidence $(x,e)$ once whenever
		$(e\setminus\{x\})\cap W\ne\varnothing$.
		Our count assigns multiplicity
		$|(e\setminus\{x\})\cap W|\le k-1$ to that incidence.
		Apply the cited lemma with parameter $\rho/(k-1)$
		and choose its corresponding constant $\lambda$.
		Multiplying its bound by $k-1$ gives \textup{(ii)}.
	\end{proof}

	\subsection{Rotation--extension for hypergraphs}
	We adapt P\'osa's rotation--extension method~\cite{Posa1962} to Berge paths.
	Let $G$ be a $k$-graph on $V$. A \emph{Berge path on $V$} is a sequence
	$
	P=v_0e_1v_1\cdots e_sv_s
	$
	with distinct \emph{core vertices} $v_0,\ldots,v_s\in V$, distinct edges
	$e_1,\ldots,e_s\in\binom Vk$, and $\{v_{i-1},v_i\}\subseteq e_i$
	for $1\le i\le s$.
	Its length is $\ell(P)=s$. Write $E(P)=\{e_1,\ldots,e_s\}$ and
	$V_{\rm core}(P)=\{v_0,\ldots,v_s\}$.
	Paths of length zero are allowed.
	A Berge path on $V$ need not be spanning, and its hyperedges
	need not belong to $G$.
	A path \emph{in $G$} must have all its hyperedges in $G$.
	The union $G\cup P$ has vertex set $V$ and edge set
	$E(G)\cup E(P)$.
	Throughout, path vertices mean core vertices.
	Berge cycles have at least two core vertices.

	Suppose that $e'\in E(G)\setminus E(P)$ contains $v_0$ and a vertex
	$v'\in V\setminus V_{\rm core}(P)$. Then
	$v'e'v_0e_1v_1\cdots e_sv_s$ is a Berge path of length $s+1$.
	We call this an \emph{extension at $v_0$}. This operation is allowed
	even if $v'$ occurs in a hyperedge of $P$ as a noncore vertex.
	
	We next describe rotations in $G$, keeping $z=v_s$ fixed.
	Our ordinary rotation is the reverse-order form of Chen and
	Zhang~\cite[Section~3, Eq.~(3.1)]{ChenZhang}.
	If an \emph{unused} edge $e''\in E(G)\setminus E(P)$ contains $v_0$ and
	$v_i$, where $1\le i\le s$, then
	\[
	v_{i-1}e_{i-1}\cdots v_1e_1v_0e''v_i
	e_{i+1}\cdots e_sv_s
	\]
	is a Berge path of length $s$, with \emph{moving endpoint} $v_{i-1}$ and
	\emph{fixed endpoint} $z$. This is an \emph{ordinary rotation}, with \emph{pivot}
	$v_i$ and \emph{broken edge} $e_i$. We also allow $i=1$, which changes only
	the first representative edge and keeps the moving endpoint $v_0$.
	Rotations from the other endpoint are defined by reversing the path and
	applying the same rule; with pivot $v_i$ they break $e_{i+1}$.
	
	A \emph{used-edge rotation} uses an edge of $G$ already in the path.
	It is the reverse-order form of Chen and Zhang's old-edge
	rotation~\cite[Section~3, Eq.~(3.2)]{ChenZhang}.  
	Write the current path from its moving endpoint $x$ as
	$P_x=x f_1u_1\cdots f_su_s$, $u_s=z.$
	If $2\le i\le s$, $f_i\in E(G)$, and $x\in f_i$, then
	\[
	u_{i-1}f_{i-1}\cdots u_1f_1x f_i u_i
	f_{i+1}u_{i+1}\cdots f_su_s
	\]
	is a path of the same length, with moving endpoint $u_{i-1}$ and fixed
	endpoint $z$. For $i=2$, the reversed prefix is $u_1f_1x$.
	Its core vertex set and edge set are unchanged. Only the pair assigned
	to $f_i$ changes, from $\{u_{i-1},u_i\}$ to $\{x,u_i\}$;
	see Salia~\cite[Section~2, Permutation~1]{Salia}.
	Both rotation types are allowed. Edges outside $G$ are not eligible
	for used-edge rotations.
	
	An original core connection $v_{j-1}e_jv_j$ is \emph{destroyed} when
	its assigned unordered pair changes or its edge is removed. Reversing
	its orientation does not destroy it. At its first destruction, one of
	$v_{j-1},v_j$ becomes the moving endpoint. Each rotation first destroys
	at most one original connection.
	Any original segment whose connections remain intact stays contiguous,
	in its original or reverse orientation, throughout each rotation sequence.
	
	An \emph{endpoint family} $\mathcal P_X=\{P_x:x\in X\}$ specifies a
	representative path and a rotation sequence in $G$ from the common path $P$
	to each moving endpoint $x$, always keeping $z$ fixed.
	Let $R_{G,z}(\mathcal P_X)$ be $X$ together with all moving endpoints
	obtained by one further allowed rotation of a selected representative.
	It depends on the representatives, not just on $X$.
	The common path, fixed endpoint, and rotation sequences are part of the family.
	For $u\ne x$ on the core of $P_x$, let $\pi_x(u)$ be its core neighbour
	towards $x$ along $P_x$.
	
	For a core vertex $x$ of a path $Q$, define its \emph{available neighbourhood}
	relative to $G$ by
	\[
	N_G(x\mid Q):=
	\{u\in V\setminus\{x\}:\exists e\in E(G)\setminus E(Q),
	\ \{x,u\}\subseteq e\}.
	\]
	This is a neighbourhood in all of $V$. If $x$ is an endpoint with no
	extension through an unused edge of $G$, then
	$N_G(x\mid Q)\subseteq V_{\rm core}(Q)$.
	For an endpoint family, set
	\[
	N_G(X\mid\mathcal P_X):=
	\left(\bigcup_{x\in X}N_G(x\mid P_x)\right)\setminus X.
	\]
	The family is \emph{ancestral} if every moving endpoint in every retained
	rotation sequence belongs to $X$, including its initial and final endpoints.
	This concerns endpoint labels: intermediate paths need not equal the
	selected representatives of those labels. Rotation sequences need not be shortest.
	
	An ancestral family can be enlarged and truncated as follows.
	If $|X|\le b\le|R_{G,z}(\mathcal P_X)|$, choose
	$X\subseteq A\subseteq R_{G,z}(\mathcal P_X)$ with $|A|=b$.
	Keep every old representative and rotation sequence. For each new endpoint,
	append one rotation to the rotation sequence of a parent in $X$.
	The resulting family is ancestral. Its rotation sequences have length at most
	$t+1$ if the old rotation sequences have length at most $t$.
	
	For the following bound, assume the family is ancestral and no selected
	$P_x$ extends at $x$ through an unused edge of $G$. This assumption holds for longest paths with fixed endpoint $z$ in
	$G\cup P$, and for longest paths in a connected component of $G$.
	Then
	\begin{equation}\label{eq:ancestral-unused-neighbourhood}
	|N_G(X\mid\mathcal P_X)|\le2|R_{G,z}(\mathcal P_X)|+3|X|.
	\end{equation}
	First, nonextension puts this neighbourhood inside $V_{\rm core}(P)$.
	Let $X^-$ and $X^+$ be the predecessor and successor sets of $X$ in
	the original core order. Let $w$ be a vertex of this neighbourhood
	outside $X\cup X^-\cup X^+$. Neither original connection next to $w$
	was destroyed in a retained
	rotation sequence. Its first destruction would put $w$ or an original neighbour
	of $w$ in $X$, a contradiction.
	An available edge from a selected endpoint to $w$ therefore gives an
	ordinary rotation to an original neighbour of $w$.
	Each resulting endpoint has at most two such original neighbours.
	These vertices contribute at most $2|R_{G,z}(\mathcal P_X)|$;
	the remaining vertices contribute at most $3|X|$.

	\begin{lem}\label{lem-expansion-adjc}
		Fix $k\ge3$ and $\varepsilon>0$.
		There are constants $C,c>0$ depending only on $k,\varepsilon$.
		For every fixed $C_0\ge C$, take $p=C_0\log n/n^{k-1}$.
		Asymptotically almost surely, $H=H^{(k)}(n,p)$ has the following properties
		for every $H'\in\mathcal H_D^\varepsilon(H)$.
		Set $G=H-H'$. Let $P$ be a Berge path on $V(H)$ with fixed endpoint
		$z$, and let $\mathcal P_X$ be a specified ancestral endpoint family
		obtained by rotations in $G$. Assume that every representative $P_x$
		admits no extension at its moving endpoint $x$ using an unused edge of $G$.
		Then
		\begin{enumerate}[label=(\roman*)]
			\item if $1\le |X|\le \lfloor(\log n)^{-1/4}/(n^{k-2}p)\rfloor$, then $|R_{G,z}(\mathcal P_X)|\ge cD|X|$;
			\item if $|X|=\lceil n/\sqrt{\log n}\rceil$, then $|R_{G,z}(\mathcal P_X)|\ge\varepsilon n/10$;
			\item the hypergraph $G$ is connected.
		\end{enumerate}
	\end{lem}
	
	\begin{proof}Let $m_1= \lfloor(\log n)^{-1/4}/(n^{k-2}p)\rfloor$, $m_2=\lceil n/\sqrt{\log n}\rceil$.
		If $\varepsilon\ge1/(4k)$, the admissible deletion family is empty for
		all sufficiently large $n$, so assume $0<\varepsilon<1/(4k)$.
		Combining \cref{lem-adjac}(i)  with the deletion condition at index $1$ gives
		\begin{equation}\label{eq:expansion-min-degree}
		\delta(G)\ge
		\left(n^{-(k-1)}+(4k-1)\varepsilon\right)D
		\ge3k\varepsilon D.
		\end{equation}
		
		We need two uniform properties of $H$. For sets $A,B\subseteq[n]$,
		apply \cref{lem-adjac} with parameter
		$2k\varepsilon/(k-1)$, and let $\lambda>0$ be its constant.
		Choose $0<L\le\min\{\lambda,\varepsilon/[16(k-1)]\}$ and let $c=L/6$. 
		The first property follows from \cref{lem-adjac}(ii): 
		\begin{equation}\label{eq:small-pair-incidences}
		J_H(A,B)\le k\varepsilon D|A|
		\quad\text{if}\quad1\le|A|\le m_1,\quad |B|\le LD|A|.
		\end{equation}
		Indeed,  $D\le2D_0$ for large $n$.
		Thus $|B|\le LD|A|\le2\lambda D_0|A|$.
		Also $m_1\le\lfloor n/((\log n)^{1/4}D_0)\rfloor$, since
		$D_0\le n^{k-1}p$.
		All set-size hypotheses of the cited lemma therefore hold.
		With the chosen parameter, its uniform bound gives
		$J_H(A,B)<k\varepsilon D_0|A|\le k\varepsilon D|A|$.
		This proves~\eqref{eq:small-pair-incidences} a.a.s.
		
		The second uniform property is
		\begin{equation}\label{eq:medium-pair-incidences}
		J_H(A,B)\le k\varepsilon D|A|
		\quad\text{if}\quad
		|A|=m_2,\qquad |B|\le\varepsilon n/2.
		\end{equation}
		For fixed $A,B$, write $J_H(A,B)=\sum_e w_e\xi_e$.
		The edge indicators $\xi_e$ are independent, and $0\le w_e\le k(k-1)$.
		Its mean satisfies
		\[
		\mu:=\mathbb E(J_H(A,B))
		\le|A||B|\binom{n-2}{k-2}p
		\le\frac{(k-1)\varepsilon}{2}D|A|
		<\frac{k\varepsilon}{2}D|A|.
		\]
		Convexity gives
		$ \mathbb E( e^{\theta J_H(A,B)})
		\le\exp\left(\frac{e^{\theta k(k-1)}-1}{k(k-1)}\mu\right).$
		Take $\theta=(\log2)/(k(k-1))$. Then
		\[
		\Pr(J_H(A,B)>k\varepsilon D|A|)
		\le\exp\left(-\frac{\log2-1/2}{k(k-1)}k\varepsilon D|A|\right)
		=\exp\{-\Omega_{k,\varepsilon}(n\sqrt{\log n})\}.
		\]
		A union bound over at most $4^n$ ordered pairs proves
		\eqref{eq:medium-pair-incidences}.
		
		Condition on these events. Fix $G$, $P$ and $\mathcal P_X$ as in the lemma.
		Write  $R=R_{G,z}(\mathcal P_X)$ and
		$S=X\cup N_G(X\mid\mathcal P_X)$. By
		\eqref{eq:ancestral-unused-neighbourhood},
		\begin{equation}\label{eq:small-rotation-cover}
		|S|\le2|R|+4|X|\le6|R|.
		\end{equation}

		Let $U$ and $B$ count the incidences $(x,e)$, $x\in X\cap e$,
		$e\in E(G)$, for which $e$ is unused or used in $P_x$, respectively.
		For an unused incidence, nonextension at $x$ puts every vertex of $e$
		on the core of $P_x$. Every vertex in $e\setminus X$ is therefore
		in the available neighbourhood, so $e\subseteq S$.
		It follows that $(k-1)U\le J_H(X,S)$.
		For a used incidence, the first edge of $P_x$ accounts for at most one
		exception per $x$. Every other such edge gives a used-edge rotation to
		a vertex $y\in R\setminus\{x\}$ with $\{x,y\}\subseteq e$.
		Thus $B\le J_H(X,R)+|X|$. Together with
		\eqref{eq:expansion-min-degree}, these bounds give
		\begin{equation}\label{eq:rotation-degree-incidences}
		3k\varepsilon D|X|\le\sum_{x\in X}d_G(x)
		=U+B\le\frac{J_H(X,S)}{k-1}+J_H(X,R)+|X|.
		\end{equation}
		
		For (i), suppose $|X|\le m_1$ but $|R|<cD|X|$.
		By~\eqref{eq:small-rotation-cover} and $6c=L$, both $|S|$ and $|R|$
		are at most $LD|X|$. Apply~\eqref{eq:small-pair-incidences} to
		$(X,S)$ and $(X,R)$. The right side of
		\eqref{eq:rotation-degree-incidences} is at most
		\[
		\left(1+\frac1{k-1}\right)k\varepsilon D|X|+|X|
		\le2k\varepsilon D|X|+|X|<3k\varepsilon D|X|
		\]
		for large $n$, a contradiction.
		
		For (ii), let $|X|=m_2$ and suppose $|R|<\varepsilon n/10$.
		Then $|S|\le2|R|+4|X|<\varepsilon n/2$ for large $n$.
		Apply~\eqref{eq:medium-pair-incidences} to $(X,S)$ and $(X,R)$.
		The same degree-incidence bound gives a contradiction.
		
		For (iii), first we show that every component has at least $\varepsilon n/10$
		vertices. Start with a longest path in a component and its initial endpoint.
		Every rotation stays in that component. An extension of any resulting path
		would contradict its maximal length there. Thus (i) and (ii) apply.
		While the ancestral family has size at most $m_1$, (i) multiplies its size
		by at least $cD$. Retain exactly $m_1$ endpoints,
		keeping the preceding family and its rotation sequences.
		One further step gives at least $cDm_1>m_2$ endpoints for large $n$.
		Since $m_1<m_2$, retain $m_2$ endpoints in the same way and apply (ii).
		This gives at least $\varepsilon n/10$ endpoints in the component.
		The construction uses $O(\log n/\log\log n)$ rotations.
		
		Suppose $G$ is disconnected. Choose a smallest component $W$, so that
		$\varepsilon n/10\le|W|\le n/2$. We need the following uniform estimate.
		\begin{clm}
			Asymptotically almost surely, every $A\subseteq W\subseteq[n]$ with
			$|A|=\lfloor\varepsilon n/40\rfloor$ and
			$\varepsilon n/10\le|W|\le n/2$ satisfies
			\[
			\sum_{v\in A}d_{H[W]}(v)
			\le(1+\varepsilon)|A|\binom{|W|-1}{k-1}p.
			\]
		\end{clm}
		\begin{proof}
			For fixed $A,W$, the left side is
			$\sum_{e\in\binom Wk}|e\cap A|\xi_e$.
			The summands are independent and bounded by $k$. Their total mean is
			$
			|A|\binom{|W|-1}{k-1}p=\Omega_{k,\varepsilon}(n\log n).
			$
			The bounded-summand Chernoff inequality gives failure probability
			$\exp\{-\Omega_{k,\varepsilon}(n\log n)\}$.
			A union bound over at most $3^n$ pairs proves the claim.
		\end{proof}
		
		Condition on this event. Let $m=|W|-\lfloor\varepsilon n/40\rfloor$,
		and order the deleted degrees as in Definition~\ref{def-resi}.
		For large $n$, $1\le m<n/2$. At most $m-1$ vertices of $W$ precede
		index $m$. Choose $A\subseteq W\cap\{v_m,\ldots,v_n\}$ with
		$|A|=\lfloor\varepsilon n/40\rfloor$.
		The deletion condition and~\cref{lem-adjac}(i)  give
		\[
		d_G(u)\ge
		\left(\left(\frac mn\right)^{k-1}+(4k-1)\varepsilon\right)D
		\qquad(u\in A).
		\]
		Since $\alpha=|W|/n\le1/2$, we have $m/n=|W|/n-|A|/n\in[0,1/2]$.
		On this interval, the derivative satisfies $(k-1)u^{k-2}\le(k-1)/2^{k-2}\le1$ for $k\ge3$.
		The mean value theorem therefore gives
		$
		\left(\frac mn\right)^{k-1}
		\ge(|W|/n)^{k-1}-\frac{|A|}{n}
		\ge(|W|/n)^{k-1}-\frac{\varepsilon}{40}.
		$
		Consequently,
		\begin{equation}\label{eq:component-lower-repaired}
		\frac1{|A|D}\sum_{u\in A}d_G(u)
		\ge(|W|/n)^{k-1}+(4k-2)\varepsilon.
		\end{equation}
		Every edge of $G$ meeting $W$ is contained in $W$. The claim gives
		\begin{align*}
		\frac1{|A|D}\sum_{u\in A}d_G(u)
		&\le(1+\varepsilon)
		\frac{\binom{|W|-1}{k-1}}{\binom n{k-1}}=(1+\varepsilon)\bigl((|W|/n)^{k-1}+O_k(1/n)\bigr)\\
		&\le(|W|/n)^{k-1}+\frac{\varepsilon}{2^{k-1}}
		+O_{k,\varepsilon}(1/n).
		\end{align*}
		Here $(|W|/n)^{k-1}\le2^{-(k-1)}$.
		Since $4k-2>2^{-(k-1)}$, this contradicts
		\eqref{eq:component-lower-repaired} for large $n$.
		Thus $G$ is connected.
	\end{proof}
	
	\section{Proof of Theorem~\ref{thm-min}}
	
	We prove Theorem~\ref{thm-min} in two steps.
	After each admissible deletion, we find a sparse spanning subhypergraph
	with the rotation--extension property.
	We then prove the required complementing property.
	Proposition~\ref{thm-comple-hami} gives Berge Hamiltonicity.

	\subsection{Subgraphs with the RE property}
	
	\begin{defi}[RE-property {\cite[Definition~2.1]{ChenZhang}}]\label{def-re}
		Let $\alpha>0$, and let $H$ be a connected $k$-graph on a vertex set
		$V$ of size $n$. 
		We say that $H$ satisfies $RE(\alpha)$ if, for every Berge path $P$
		on $V$, at least one of the following holds:
		\begin{enumerate}[label=(\roman*)]
			\item There exists a Berge path in $H\cup P$ of length longer than
			$\ell(P)$.
			\item There exists a set $S_P\subseteq V$ with $|S_P|\ge\alpha n$.
			For every $v\in S_P$, there exists a set
			$T_v\subseteq V\setminus\{v\}$ with $|T_v|\ge\alpha n$.
			For every $w\in T_v$, there is a Berge path
			$Q_{v,w}\subseteq H\cup P$ with endpoints $v,w$ and
			$\ell(Q_{v,w})=\ell(P)$.
		\end{enumerate}
	\end{defi}
	
	For each Berge path $P$, the RE property therefore provides either a
	longer path or large sets of endpoint pairs joined by paths of length $\ell(P)$. Theorem~\ref{thm-edge-re} gives a sparse spanning subhypergraph with property $RE(1/2+\varepsilon)$. We prove it using the endpoint-growth estimates below.

	\begin{lem}\label{lem:degree-endpoint-growth}
		Fix $k\ge3$, $\varepsilon>0$, and a constant $A>0$.
		Let $t_n:=\lceil A\log n/\log\log n\rceil$.
		There exists $C=C(k,\varepsilon,A)>0$ with the following property.
		For $C\log n/n^{k-1}\le p\le1$, a.a.s. $H=H^{(k)}(n,p)$ satisfies
		the result below  for all the following choices.
		Let $H'\subseteq H$ and  $G:=H-H'$.
		Let $P'$ be any Berge path on $[n]$, and let $P$ be a longest Berge
		path in $G\cup P'$. Fix an endpoint $z$ of $P$.
		Let $X$ be a set of at least $\varepsilon n/20$ distinct moving endpoints.
		For each $x\in X$, choose a path $Q_x$ obtained from $P$ by at most
		$t_n$ allowed rotations in $G$, keeping $z$ fixed.
		Write $\mathcal Q_X:=\{Q_x:x\in X\}$ and use
		$R_{G,z}(\mathcal Q_X)$ as defined above.
		This family need not be ancestral.
		If $a\in[0,1]$ and
		$
		d_{H'}(x)\le (1-a-4k\varepsilon)D$ for every $x\in X,
		$
		then
		\[
		\left(\frac{|R_{G,z}(\mathcal Q_X)|}{n}\right)^{k-1}
		\ge a+k\varepsilon-o(1).
		\]
	\end{lem}
	
	\begin{proof}
		We first establish a uniform probabilistic estimate using a deterministic
		family of encodings. This step is needed because the paths $Q_x$ and
		the endpoint set $R_{G,z}(\mathcal Q_X)$ may depend on the random hypergraph.
		Put $L:=\lceil32t_n/\varepsilon\rceil$.
		Write the original path as $P=v_0e_1v_1\cdots e_\ell v_\ell$,
		oriented so that $z=v_\ell$.
		Since $\ell+1\ge |X|\ge\varepsilon n/20$, we have $L<\ell+1$ for all sufficiently
		large $n$. Use a fixed rule to partition the ordered core vertex set
		into $L$ consecutive nonempty blocks, each of size at most
		$\lceil (\ell+1)/L\rceil$.
		
		For a chosen rotation sequence leading to $Q_x$, call a block bad if
		an original core adjacency pair internal to that block has been destroyed
		at some step of the sequence. An original core adjacency pair here means
		the occurrence $v_{i-1}e_iv_i$, without an orientation. Ordinary
		rotations remove one such current adjacency, while a used-edge rotation
		changes the core pair represented by one current edge. Thus at most $t_n$ original adjacency pairs have ever been destroyed.
		At most $t_n$ blocks are therefore bad. All original internal adjacencies of every
		other block remain intact throughout the sequence. Such a block therefore
		appears as a contiguous segment of $Q_x$, either in its original order
		or in the reverse order.
		
		Let $E_x$ consist of all vertices of the bad blocks, the first and last
		vertices of every block, and $z$. Then, uniformly in all the choices,
		\[
		|E_x|\le t_n\lceil (\ell+1)/L\rceil+2L+1
		\le \varepsilon n/32+t_n+2L+1\le\varepsilon n/8
		\]
		for sufficiently large $n$. On the internal vertices of each nonbad
		block, the neighbour towards $x$ along $Q_x$ is either the successor
		or the predecessor in the original core order, with one direction for
		the entire block. Denote this map by $f_x$ on those internal vertices.
		In particular, the moving endpoint $x$ cannot be an internal vertex of
		a nonbad block.
		
		For any set $B\subseteq[n]$, define
		\[
		\widetilde T_x:=
		\left(E_x\cup
		\{u\notin E_x:\ u\in V_{\rm core}(P),\ f_x(u)\in B\}\right)
		\setminus\{x\}.
		\]
		The map $f_x$ is injective on the interior of each nonbad block and
		maps that interior into the same block. Its images in distinct blocks
		are disjoint. Hence
		\begin{align}  \label{*}
		|\widetilde T_x|\le |B|+\varepsilon n/8.               
		\end{align}
		If $B=R_{G,z}(\mathcal Q_X)$ and $\pi_x(u)$ denotes the neighbour of $u$
		towards $x$ along $Q_x$, then the set
		\[
		T_x:=\{z\}\cup
		\{u\in V_{\rm core}(P)\setminus\{x,z\}:\ \pi_x(u)\in B\}
		\]
		satisfies $T_x\subseteq\widetilde T_x$.
		
		An \emph{encoding} records the ordered core vertices of the original path, the sets $X$ and $B$, and, for each endpoint in $X$, whether each block is bad, forward or reverse. It determines the sets $\widetilde T_x$ used in the probabilistic estimate. We only consider encodings with $|X|\ge\varepsilon n/20$
		and at most $t_n$ bad blocks for each $x\in X$.
	  There are at most
		\begin{align}\label{**}
		\left(\sum_{\ell=0}^{n-1} n^{\ell+1}\right)4^n3^{nL}
		\le n^{n+1}4^n3^{nL}
		=\exp\bigl((1+o(1))n\log n\bigr)               
		\end{align}
		such encodings, since $L=O_{\varepsilon,A}(\log n/\log\log n)$.
		The construction and bound \eqref{*} make sense for each of these fixed
		encodings independently of the random hypergraph.
		
		Fix one encoding  and let $\beta:=|B|/n$. For each $x\in X$ define
		\[
		\widetilde{\mathcal D}_x:=
		\{e\in\tbinom{[n]}k:\ x\in e,
		\ e\setminus\{x\}\not\subseteq\widetilde T_x\},
		\qquad
		\widetilde Z:=\sum_{x\in X}|E(H)\cap\widetilde{\mathcal D}_x|.
		\]
		Since $x\notin\widetilde T_x$, we have exactly
		$
		|\widetilde{\mathcal D}_x|
		=\binom{n-1}{k-1}-\binom{|\widetilde T_x|}{k-1}.
		$
		Using \eqref{*} and the fact that $u\mapsto u^{k-1}$ has derivative at most
		$k-1$ on $[0,1]$, we obtain
		\[
		\begin{split}
		\mu:=\mathbb E(\widetilde Z)
		&\ge\left(1-\min\{1,\beta+\varepsilon/8\}^{k-1}-O_k(n^{-1})\right)|X|D\\
		&\ge\left(1-\beta^{k-1}-(k-1)\varepsilon/8-O_k(n^{-1})\right)|X|D,
		\end{split}
		\]
		whereas $\mu\le |X|D$.
		
		For each fixed encoding, let $\widetilde Z:=\sum_e w_e\xi_e$.
		The edge indicators are independent, and the weights are deterministic
		integers with $0\le w_e\le k$.
		If $1-\beta^{k-1}-2k\varepsilon\le0$, the required lower bound
		$\widetilde Z\ge(1-\beta^{k-1}-2k\varepsilon)|X|D$ is automatic.
		Otherwise the expectation bound gives
		$\mu-(1-\beta^{k-1}-2k\varepsilon)|X|D\ge k\varepsilon|X|D>0$
		for sufficiently large $n$. 
		Apply the weighted lower-tail bound in~\cite[Theorem~2.6]{ChenZhang}
		with additive deviation $k\varepsilon|X|D$.
		Since $\mu\le|X|D$, it gives
		\[
		\Pr\left(\widetilde Z<
		(1-\beta^{k-1}-2k\varepsilon)|X|D\right)
		\le\exp\left(-\frac{k\varepsilon^2}{2}|X|D\right).
		\]
		Since $|X|\ge\varepsilon n/20$ and
		$D\ge (1-o(1))C\log n/(k-1)!$, taking $C$ sufficiently large
		makes the union bound over \eqref{**} tend to zero. Therefore, for every
		encoding,
		\begin{align}\label{***}
		\widetilde Z\ge
		(1-\beta^{k-1}-2k\varepsilon)|X|D.         \end{align}
		
		Now fix arbitrary choices in the lemma, possibly depending on $H$.
		Use their encoding with $B=R_{G,z}(\mathcal Q_X)$.
		Suppose that
		$e\in E(G)\cap\widetilde{\mathcal D}_x$ and $e\notin E(Q_x)$.
		There exists $u\in e\setminus(\{x\}\cup\widetilde T_x)$. If $u$
		is outside $V_{\rm core}(Q_x)=V_{\rm core}(P)$, adding $u$ through
		$e$ extends $Q_x$ in $G\cup P'$, contrary to the maximality of $P$.
		Otherwise $u$ is a core vertex. It is not $z$, since
		$z\in\widetilde T_x$. It is also not the core neighbour of $x$ along
		$Q_x$. That neighbour maps to $x\in X\subseteq B$ under $\pi_x$,
		so it belongs to $T_x\subseteq\widetilde T_x$.
		Thus an ordinary unused-edge rotation with pivot $u$ is valid. Its new endpoint is
		$\pi_x(u)\notin B$, since $u\notin T_x$. This contradicts the
		definition of $B=R_{G,z}(\mathcal Q_X)$.
		
		Every pair counted by $\widetilde Z$ must therefore have its edge
		in $H'$ or in the corresponding $Q_x$. Hence
		\[
		\widetilde Z\le
		\sum_{x\in X}d_{H'}(x)
		+\sum_{x\in X}|\{e\in E(Q_x):x\in e\}|.
		\]
		The common path contributes at most $kn$ incidences.
		Each $Q_x$ has at most $t_n$ additional representatives, since each
		rotation adds at most one and used-edge rotations add none. Since $|X|\ge\varepsilon n/20$ and $t_n=o(D)$,
		\[
		\sum_{x\in X}|\{e\in E(Q_x):x\in e\}|
		\le kn+t_n|X|=o(|X|D).
		\]
		Combining the degree hypothesis with \eqref{***} now gives
		\[
		(1-\beta^{k-1}-2k\varepsilon)|X|D
		\le (1-a-4k\varepsilon+o(1))|X|D.
		\]
		Therefore $\beta^{k-1}\ge a+2k\varepsilon-o(1)
		\ge a+k\varepsilon-o(1)$, as required.
	\end{proof}
	
	\begin{lem}\label{lem-re}
		Fix $k\ge3$ and $\varepsilon>0$.
		There exists a constant $C_*>0$ with the following property.
		For every fixed $C_0\ge C_*$ and $p=C_0\log n/n^{k-1}$, a.a.s.
		$H=H^{(k)}(n,p)$ satisfies the result below.
		For every $H'\in\mathcal H_D^\varepsilon(H)$, the hypergraph
		$H''=H-H'$ satisfies $RE(1/2+\varepsilon)$.
	\end{lem}
	\begin{proof}
		If $\varepsilon\ge1/(4k)$, the deletion class is empty. We therefore
		assume $0<\varepsilon<1/(4k)$. Fix 
		$t_n=\lceil4\log n/\log\log n\rceil$ before choosing $C_0$.
		Choose $C_0$ sufficiently large for all
		the uniform estimates below. Condition first on Lemma~\ref{lem-expansion-adjc}.
		In particular $H''$ is connected.
		
		Let $P'$ be an arbitrary Berge path on $V(H)$. If $H''\cup P'$ contains
		a longer Berge path, alternative (i) of Definition~\ref{def-re} holds.
		Otherwise let $P=P'=v_0e_1\cdots e_sv_s$, and keep $v_s$ fixed.
		All rotations below use $H''$ and the two operations defined above.
		
		Start with the endpoint  $X_0=\{v_0\}$. At each step retain the
		old representative paths and give each newly selected endpoint a path
		obtained by one rotation from an old representative. While its size is at most
		$m_1=\lfloor(\log n)^{-1/4}(n^{k-2}p)^{-1}\rfloor$,
		Lemma~\ref{lem-expansion-adjc}(i) multiplies its size by at least $cD$.
		If a step first exceeds $m_1$, retain exactly $m_1$ endpoints, including
		all endpoints from the preceding step.  Thus after
		$O(\log n/\log\log n)$ rotations, we have a family
		of size $m_1$ with the required rotation sequences. One more application
		gives at least $cDm_1\gg n/\sqrt{\log n}$ endpoints. Retain exactly
		$\lceil n/\sqrt{\log n}\rceil$ of them.
		Applying Lemma~\ref{lem-expansion-adjc}(ii) now gives an endpoint family
		of size at least $\varepsilon n/10$.
		
		We next grow a linear-sized family. Note that each
		representative uses at most $t_n$ rotations from the common path $P$.
		Indeed $cD\ge\sqrt{\log n}$ for sufficiently large $n$, so the initial
		construction uses at most $2\log n/\log\log n+O(1)$ rotations.

		Let $U$ be the current endpoint set, and write $\alpha=|U|/n$.
		We have $\alpha\ge\varepsilon/10$. Let
		$
		m=\min\{|U|,\lfloor n/2\rfloor-1\}
		-\lfloor\varepsilon n/20\rfloor.
		$
		Order the deleted degrees as in Definition~\ref{def-resi}. At most
		$m-1$ vertices of $U$ have indices below $m$. We can therefore choose
		$X\subseteq U\cap\{v_m,\ldots,v_n\}$ of size
		$\lceil\varepsilon n/20\rceil$, adjusting $m$ by one if necessary.
		For every $x\in X$, the nonlinear deletion condition gives
		\[
		d_{H'}(x)\le
		\left(1-\left(\frac mn\right)^{k-1}-4k\varepsilon\right)D.
		\]
		Retain the chosen representative for each $x\in X$ and denote their
		indexed family by $\mathcal P_X$. Apply
		Lemma~\ref{lem:degree-endpoint-growth} with $G=H''$, $z=v_s$,
		$Q_x=P_x$, and $a=(m/n)^{k-1}$. Let the next family have endpoint set $R_{H'',v_s}(\mathcal P_X)$.
		Write $\beta=|R_{H'',v_s}(\mathcal P_X)|/n$.
		The lemma gives $\beta^{k-1}\ge(m/n)^{k-1}+k\varepsilon-o(1)$.
		For sufficiently large $n$, we have
		$0\le m/n\le\min\{\alpha,1/2\}\le1/2$.
		The definition of $m$, including its possible adjustment by one, gives
		$\min\{\alpha,1/2\}-m/n=\varepsilon/20+O(1/n)$.
		Since $(k-1)t^{k-2}\le1$ on $[0,1/2]$, the mean value theorem yields
		\[
		\left(\frac mn\right)^{k-1}
		\ge\min\{\alpha,1/2\}^{k-1}-\frac{\varepsilon}{20}-O_k(1/n).
		\]
		Thus $\beta^{k-1}\ge\min\{\alpha,1/2\}^{k-1}
		+(k-1/20)\varepsilon-o(1)$.
		For all sufficiently large $n$, this implies
		\begin{equation}\label{eq-ration}
		\beta^{k-1}\ge
		\min\{\alpha,1/2\}^{\,k-1}
		+\left(k-\frac12\right)\varepsilon.
		\end{equation}
		If $\varepsilon/10\le\alpha\le1/2$, this implies
		$\beta^{k-1}>\alpha^{k-1}$ and hence $\beta>\alpha$. Unless
		$\beta\ge1/2+\varepsilon$, the mean
		value theorem on $[\alpha,\beta]\subseteq[0,1]$ therefore gives
		$
		(k-1)(\beta-\alpha)\ge
		\beta^{k-1}-\alpha^{k-1}\ge
		\left(k-\frac12\right)\varepsilon.
		$
		Thus each step increases the family size by more than $\varepsilon n$
		until it reaches $n/2$. One final application of \eqref{eq-ration}
		gives
		$
		\beta^{k-1}\ge 2^{-(k-1)}+
		\left(k-\frac12\right)\varepsilon
		\ge (1/2+\varepsilon)^{k-1},
		$
		where the last inequality follows from 
		$(k-1)u^{k-2}\le k-1$ for $0\le u\le1$.
		The total number of additional rotations is $O_\varepsilon(1)$,
		as allowed in the choice of $t_n$.
		
		We have therefore obtained $S_P$ of size at least
		$(1/2+\varepsilon)n$. Each $v\in S_P$ is joined to $v_s$ by a path
		of the same length as $P$. Fix such a $v$ and take its representative
		path as a new common path. Repeat the argument keeping $v$ fixed.
		All uniform estimates apply to this new path as well.
		This gives $T_v$ of size at least $(1/2+\varepsilon)n$, with an
		equal-length path joining $v$ to every $w\in T_v$ in $H''\cup P'$.
		Together with the connectivity established in
		Lemma~\ref{lem-expansion-adjc}(iii), this is $RE(1/2+\varepsilon)$.
	\end{proof}
	
	\begin{thm}\label{thm-edge-re}
		Fix an integer $k\ge3$, $\varepsilon>0$, and a constant
		$0<\delta<1$. There exists $C=C(k,\varepsilon,\delta)>0$ with the
		following property. If $C\log n/n^{k-1}\le p\le1$, then a.a.s.
		$H=H^{(k)}(n,p)$ satisfies the result below.
		For every $H'\in\mathcal H_D^{2\varepsilon}(H)$, the hypergraph
		$H-H'$ contains a spanning subhypergraph on $[n]$ with at most
		$\delta n^kp$ edges satisfying $RE(1/2+\varepsilon)$.
	\end{thm}
	The proof of Theorem~\ref{thm-edge-re} is similar to the proof of Lemma 25 in
	Condon et al.~\cite{CondonEtAl}. For completeness, we include
	the details in the appendix.

	\subsection{Complementing subgraphs}       
	
	The preceding theorem supplies a sparse spanning subhypergraph after an
	admissible deletion. We next show that the remaining hypergraph complements
	every such subhypergraph. Proposition~\ref{thm-comple-hami} then yields
	Berge Hamiltonicity.
	
	\begin{defi}[{\cite[Definition~2.3]{ChenZhang}}]\label{def-comple}
		Let $\alpha>0$, and let $H_1,H_2$ be $k$-graphs on a common vertex set
		$V$ of size $n$, with $H_1$ satisfying $RE(\alpha)$. We say that $H_2$
		\emph{complements $H_1$ with respect to $RE(\alpha)$} if every Berge
		path $P$ on $V$, as in Definition~\ref{def-re}, satisfies one of the
		following:
		\begin{enumerate}[label=(\roman*)]
			\item The hypergraph $H_1\cup P$ contains a Berge path longer than $P$.
			\item There are sets $S_P$ and $T_v$, $v\in S_P$, satisfying
			condition \textup{(ii)} of Definition~\ref{def-re} for $H_1$ and $P$.
			For some $v\in S_P$ and $w\in T_v$, there is a Berge path
			$Q_{v,w}\subseteq H_1\cup P$ with endpoints $v,w$ and length $\ell(P)$,
			and an edge $e\in E(H_2)\setminus E(Q_{v,w})$ containing $v,w$.
			We call $e$ a \emph{booster for $Q_{v,w}$}. Consequently,
			$Q_{v,w}\cup\{e\}$ contains a Berge cycle of length $\ell(P)+1$
			with core vertex set $V_{\rm core}(Q_{v,w})$.
		\end{enumerate}
	\end{defi}
	
	\begin{prop}[ {\cite[Proposition~2.5]{ChenZhang}}]\label{thm-comple-hami}
		Let $n\ge k\ge3$ and $0<\alpha<1$, and let $H_1,H_2$ be $k$-graphs
		on a common set of $n$ vertices. If $H_1$ satisfies $RE(\alpha)$ and
		$H_2$ complements $H_1$ with respect to $RE(\alpha)$, then
		$H_1\cup H_2$ is Berge Hamiltonian.
	\end{prop}

	The next theorem supplies the complementing property 
	for every sufficiently sparse $RE$-subhypergraph of the random hypergraph.
	Together with Theorem~\ref{thm-edge-re} and
	Proposition~\ref{thm-comple-hami}, it completes the proof.  
	
	\begin{thm}\label{thm-comple}
		For every fixed $k\ge3$ and $0<\varepsilon<1/2$, there are
		$\delta,C>0$ such that the following holds a.a.s. for
		$H=H^{(k)}(n,p)$ when $p\ge C\log n/n^{k-1}$.
		For every $H'\in\mathcal H_D^\varepsilon(H)$, the $k$-graph $H-H'$
		complements every spanning
		$R\subseteq H$ on $[n]$ satisfying $RE(1/2+\varepsilon)$ and
		$e(R)\le\delta n^kp$.
	\end{thm}
	\begin{proof}
		It suffices to prove the assertion for
		$
		0<\eta:=\min\{\varepsilon,1/(16k)\}.
		$
		Indeed, $\mathcal H_D^\varepsilon(H)\subseteq
		\mathcal H_D^\eta(H)$, and $\operatorname{RE}(1/2+\varepsilon)$
		implies $\operatorname{RE}(1/2+\eta)$.

		For a labelled $k$-graph $R$ satisfying
		$RE(1/2+\varepsilon)$ and
		$
		e(R)\leq\delta n^kp,
		$
		and  $R\subseteq H$. All edges outside $R$ remain
		independent Bernoulli random variables with parameter $p$.
			Fix an arbitrary Berge path $P$ on the vertex set of $R$. If $R\cup P$ contains a longer path,
		the first complementing alternative holds. Otherwise, 
		from $R$ and $P$, fix sets $S_P$ and $T_v$
		satisfying $RE(1/2+\varepsilon)$. Thus
		\[
		|S_P|\geq(1/2+\eta)n,\qquad
		|T_v|\geq(1/2+\eta)n\quad(v\in S_P).
		\]
		For every $w\in T_v$ there is a path $Q_{v,w}\subseteq R\cup P$
		of length $\ell(P)$ with endpoints $v,w$.
		
		Set $s:=\lfloor\eta n/2\rfloor$. For every $S\subseteq S_P$
		with $|S|=s$, define
		\[
		\mathcal C(S):=
		\left\{
		e\in\binom{[n]}k\setminus(E(R)\cup E(P)):
		\begin{array}{l}
		e\cap S=\{v\}\text{ for some }v\in S,\\
		(e\setminus\{v\})\cap(T_v\setminus S)\ne\varnothing
		\end{array}
		\right\}.
		\]
		If $e\in E(H-H')\cap\mathcal C(S)$, it contains some
		$v\in S$ and $w\in T_v\setminus S$. Since $Q_{v,w}\subseteq R\cup P$,
		the edge $e$ is not used by $Q_{v,w}$ and is a booster. Hence if $H-H'$ does not
		complement $R$, then for every choice of $S$ as above,
		\begin{equation}\label{eq:complement-upper}
		|E(H)\cap\mathcal C(S)|
		\le |E(H')\cap\mathcal C(S)|
		\le\sum_{v\in S}d_{H'}(v).
		\end{equation}
		
		We next establish a lower bound that holds for
		all such $S$. Write $Y_v:=T_v\setminus S$.
		For a fixed $v\in S$,
		the number of $k$-sets meeting $S$ exactly in $v$ and meeting $Y_v$
		is $\binom{n-s}{k-1}-\binom{n-s-|Y_v|}{k-1}$.
		Since $|T_v|\ge(1/2+\eta)n$, we have
		$n-s-|Y_v|\le(1/2-\eta)n$.
		The sets counted for distinct $v\in S$ are disjoint. Removing the edges
		of $R$ and $P$ over $S$ therefore gives
		\[
		|\mathcal C(S)|\ge
		s\left[\binom{n-s}{k-1}
		-\binom{\lfloor(1/2-\eta)n\rfloor}{k-1}\right]-e(R)-e(P).
		\]
		Because $\mathcal C(S)$ excludes $E(R)$, the variable
		$Z_S:=|E(H)\cap\mathcal C(S)|$ is binomial with mean
		$\mathbb E (Z_S)=p|\mathcal C(S)|$.
		Using $e(P)\le n-1$ and $p\le1$, we obtain
		\[
		\frac{\mathbb{E}(Z_S)}{sD}
		\geq
		1-{1\over 2^{k-1}}-\frac{k-1}{2}\eta
		-\frac{2(k-1)!\delta}{\eta}-o(1).
		\]
		Choose $\delta\leq\eta^2/[16(k-1)!]$. For all sufficiently
		large $n$, this implies
		\begin{equation}\label{eq:complement-mean}
		\mathbb{E}(Z_S)\geq s\left( 1-{1\over 2^{k-1}}-2k\eta\right) D.
		\end{equation}
		
		Let $H'\in\mathcal H_D^\eta(H)$, and order its vertices by
		nonincreasing deleted degree as $v_1,\ldots,v_n$.
		Since
		$|S_P|\geq(1/2+\eta)n$, at least $\eta n-O(1)$ vertices of
		$S_P$ occur at positions $i\geq\lceil n/2\rceil$. We may
		therefore choose $S\subseteq S_P$ of size $s$ entirely from
		those positions. Applying Definition~\ref{def-resi} at
		$j=\lceil n/2\rceil-1<n/2$ yields
		\[
		d_{H'}(v)\le
		\left(1-\left(\frac jn\right)^{k-1}-4k\eta\right)D
		\le\left(1-2^{-(k-1)}-3k\eta\right)D
		\qquad(v\in S).
		\]
		The last inequality holds for sufficiently large $n$, since
		$(j/n)^{k-1}=2^{-(k-1)}-O_k(n^{-1})$.
		By \eqref{eq:complement-upper}, failure to complement $R$ implies
		$Z_S\le s(1-2^{-(k-1)}-3k\eta)D$ for at least one such $S$.
		By \eqref{eq:complement-mean}, Chernoff's inequality gives
		\[
		\Pr\!\left(
		Z_S\leq s\left( 1-{1\over 2^{k-1}}-3k\eta\right) D
		\,\mid R\subseteq H
		\right)
		\leq \exp(-c_{k,\eta}n^kp)
		\]
		for some $c_{k,\eta}>0$.
		
		There are at most $\exp(O_k(n\log n))$ labelled Berge paths $P$
		and $2^n$ choices of $S$. Choose $C=C(k,\eta)$ sufficiently large.
		A union bound gives that for each fixed $R$ there exists $b_{k,\eta}>0$ such that  
		\begin{equation}\label{eq:conditional-complement}
		\Pr\bigl(
		\exists H'\in\mathcal H_D^\eta(H):
		H-H'\text{ does not complement }R
		\,\bigm|\,R\subseteq H
		\bigr)
		\leq \exp(-b_{k,\eta}n^kp).
		\end{equation}

		Finally, sum \eqref{eq:conditional-complement} over all such labelled
		$R$. The failure probability is at most
		\[
		\exp(-b_{k,\eta}n^kp)
		\sum_{m=0}^{\lfloor\delta n^kp\rfloor}\binom{ \binom nk}{m}p^m.
		\]
		For $\delta$ sufficiently small, $\delta n^kp<\binom nkp$ for large $n$, and
		\[
		\sum_{m=0}^{\lfloor\delta n^kp\rfloor}\binom {\binom nk}{m} p^m
		\leq
		(\lfloor\delta n^kp\rfloor+1)\left(\frac{e\binom nkp}{\lfloor\delta n^kp\rfloor}\right)^{\lfloor\delta n^kp\rfloor}
		=
		\exp\!\left(
		\bigl[\delta\log(e/(k!\delta))+o(1)\bigr]n^kp
		\right).
		\]
		Choosing $\delta$  smaller so that
		$\delta\log(e/(k!\delta))<b_{k,\eta}/2$,
		the total failure probability is at most
		\[
		\exp\!\left(
		\bigl[-b_{k,\eta}+\delta\log(e/(k!\delta))+o(1)\bigr]n^kp
		\right)=o(1).
		\]
		This completes the proof.
	\end{proof}

	\begin{proof}[Proof of Theorem \ref{thm-min}]
		If $\varepsilon\ge1/(4k)$, Definition~\ref{def-resi} has a negative
		upper bound at $i=1$, so the deletion class is empty. Assume that
		$0<\varepsilon<1/(4k)$.
		Choose $0<\delta<1$ small enough for Theorem~\ref{thm-comple} with
		parameter $\varepsilon/2$. Then choose $C$ large enough for that theorem
		and Theorem~\ref{thm-edge-re} with parameter $\varepsilon/2$ and this
		$\delta$.
		
		Let $H=H^{(k)}(n,p)$ and condition on both theorem events.
		Fix any $H'\in\mathcal H_D^\varepsilon(H)$ and let $H''=H-H'$.
		Theorem~\ref{thm-edge-re} supplies a spanning $H^*\subseteq H''$ on
		$[n]$ with at most $\delta n^kp$ edges and property $RE(1/2+\varepsilon/2)$.
		Since $H'\in\mathcal H_D^\varepsilon(H)\subseteq
		\mathcal H_D^{\varepsilon/2}(H)$, Theorem~\ref{thm-comple} says that
		$H''$ complements $H^*$ with respect to $RE(1/2+\varepsilon/2)$.
		By Proposition~\ref{thm-comple-hami}, the union $H''=H-H'$ is Berge Hamiltonian.
	\end{proof}

	\section{Appendix}
	\begin{proof}[Proof of Theorem~\ref{thm-edge-re}]

		Fix $0<\delta<1$. Choose $C_0$ large enough for Lemma~\ref{lem-re}
		with parameter $\varepsilon$ and for the degree bounds below. Let
		$C\ge C_0/\delta$, $p_0=C_0\log n/n^{k-1}$ and
		$q=p_0/p\le\delta$.
		Conditional on $H=H^{(k)}(n,p)$, retain its edges independently with
		probability $q$, keeping vertex set $[n]$ to obtain $\widehat H$.
		Then $\widehat H\sim H^{(k)}(n,p_0)$ and
		$D':=\binom n{k-1}p_0=qD=\Theta(\log n)$.
		
		Call a $k$-graph $F$ on $[n]$ $D'$-$good$ if $F-J$ satisfies
		$RE(1/2+\varepsilon)$ for every $J\in\mathcal H_{D'}^\varepsilon(F)$. Lemma~\ref{lem-re} implies $\widehat H$ is
		$D'$-good with probability $1-o(1)$.
		Since $q\le\delta$, we have
		$\mathbb E (e(\widehat H))=\binom nkp_0\le\delta n^kp/k!$.
		Chernoff's inequality~\cite[Theorem~2.6(ii)]{ChenZhang} gives
		$e(\widehat H)\le\delta n^kp$ with probability $1-o(1)$.
		Let $\mathcal A$ be the event that both this bound and $D'$-goodness
		hold. Thus $\Pr(\mathcal A)=1-o(1)$.
		
		Fix any  $H$ and any $H'\in\mathcal H_D^{2\varepsilon}(H)$.
		We will show that
		\begin{align}\label{17}
		\widehat H\cap H'\in\mathcal H_{D'}^\varepsilon(\widehat H).
		\end{align}
		This holds with conditional probability at least $1-n^{-2}$, uniformly
		in both choices. 
		Since $H'\in\mathcal H_{D}^{2\varepsilon}(H)$, there exists an
		ordering $v_1,\ldots,v_n$ of the vertices such that
		$
		d_{H'}(v_1)\ge d_{H'}(v_2)\ge\cdots\ge d_{H'}(v_n) 
		$
		and for every $i<n/2$,
		\begin{align}\label{19}
		d_{H'}(v_i)
		\le \left( 1
		-\left(\frac{i}{n}\right)^{k-1}
		-8k\varepsilon\right) D. 
		\end{align}
		
		Let $J=\widehat H\cap H'$. For each vertex $v\in V(J)$,
		$d_J(v)\sim\operatorname{Bin}(d_{H'}(v),q)$ and
		$\mathbb E( d_J(v))=q d_{H'}(v)\le D'$.
		The additive Chernoff bound gives
		\begin{align}\label{21}
		\Pr\!\left(
		d_J(v)>q d_{H'}(v)+2k\varepsilon D'
		\right)\le \exp\{-c_{k,\varepsilon}D'\} 
		\end{align}
		for some $c_{k,\varepsilon}>0$. Since
		$D'\ge c_kC_0\log n$, choosing $C_0$ large enough makes this
		probability at most $n^{-3}$. A union bound gives
		$d_J(v)\le qd_{H'}(v)+2k\varepsilon D'$ for every vertex, with
		probability at least $1-n^{-2}$.
		For $i<n/2$ and $j\ge i$, we have $d_{H'}(v_j)\le d_{H'}(v_i)$.
		Together with \eqref{19}, we have
		\begin{align}\label{23}
		d_J(v_j)
		\le qd_{H'}(v_i)+2k\varepsilon D'
		\le\left(1-\left(\frac in\right)^{k-1}-6k\varepsilon\right)D'.
		\end{align}
		Let $w_1,\ldots,w_n$ be an ordering of the vertices satisfying
		$
		d_J(w_1)\ge\cdots\ge d_J(w_n). 
		$
		We claim that, for every $i<n/2$,
		\begin{align}\label{25}
		d_J(w_i)
		\le \left( 1-\left(\frac{i}{n}\right)^{k-1}-4k\varepsilon\right)  D'. 
		\end{align}
		By \eqref{23}, every vertex among $v_i,\ldots,v_n$ satisfies
		the bound in \eqref{25}.
		Thus at most $i-1$ vertices exceed its right-hand side.
		The $i$-th largest degree therefore satisfies \eqref{25}. Since $J\subseteq\widehat H$, the bounds in~\eqref{25} imply~\eqref{17}.

		Apply Chen and Zhang~\cite[Lemma~4.2]{ChenZhang} with $\mathcal A$
		and, for each $H'\in\mathcal H_D^{2\varepsilon}(H)$, the event in
		\eqref{17}. Its conditional failure probability is at most $n^{-2}$
		for every fixed $H,H'$. Since $\Pr(\mathcal A)=1-o(1)$, a.a.s. over
		$H$, every such $H'$ admits a conditional sample for which both events
		hold. The sample may depend on $H'$.
		
		Its $D'$-goodness implies that
		$\widehat H-J=\widehat H-H'$ satisfies $RE(1/2+\varepsilon)$.
		This is a spanning subhypergraph of $H-H'$ on $[n]$, with at most
		$e(\widehat H)\le\delta n^kp$ edges. This proves the theorem.
	\end{proof}


\begin{thebibliography}{99}
		
		\bibitem{Karp}
		R.~M. Karp, Reducibility among combinatorial problems, in: R.~E. Miller and J.~W. Thatcher (Eds.), \emph{Complexity of Computer Computations}, Plenum Press, New York, 1972, pp.~85--103.
		
		\bibitem{Dirac}
		G.~A. Dirac, Some theorems on abstract graphs, \emph{Proc. London Math. Soc.} (3) 2 (1952), 69--81.
		
		\bibitem{Posa1962}
		L. P\'osa, A theorem concerning Hamilton lines, \emph{Magyar Tud. Akad. Mat. Kutat\'o Int. K\"ozl.} 7 (1962), 225--226.
		
		\bibitem{Chvatal}
		V. Chv\'atal, On Hamilton's ideals, \emph{J. Combin. Theory Ser. B} 12 (1972), 163--168.
	
		\bibitem{SudakovVu}
		B. Sudakov and V.~H. Vu, Local resilience of graphs, \emph{Random Structures Algorithms} 33 (2008), 409--433.
		
		\bibitem{LeeSudakov}
		C. Lee and B. Sudakov, Dirac's theorem for random graphs, \emph{Random Structures Algorithms} 41 (2012), 293--305.
	
		
	\bibitem{CondonEtAl}
	P. Condon, A. Espuny D\'iaz, D. K\"uhn, D. Osthus and J. Kim, Resilient degree sequences with respect to Hamilton cycles and matchings in random graphs, \emph{Electron. J. Combin.} 26 (2019), no.~4, Paper No.~P4.54.
	
		
		\bibitem{Berge}
		C. Berge, \emph{Graphes et hypergraphes}, Dunod, Paris, 1970.
		
		\bibitem{BermondEtAl}
		J.-C. Bermond, A. Germa, M.-C. Heydemann and D. Sotteau, Hypergraphes hamiltoniens, in: \emph{Probl\`emes combinatoires et th\'eorie des graphes}, Colloq. Internat. CNRS, Vol.~260, CNRS, Paris, 1978, pp.~39--43.
		
		\bibitem{FurediKostochkaLuo}
		Z. F\"uredi, A. Kostochka and R. Luo, Berge cycles in non-uniform hypergraphs, \emph{Electron. J. Combin.} 27 (2020), no.~3, Paper No.~P3.9.
		
		\bibitem{KostochkaLuoMcCourt}
		A. Kostochka, R. Luo and G. McCourt, Dirac-type theorems for long Berge cycles in hypergraphs, \emph{J. Combin. Theory Ser. B} 168 (2024), 159--191.
		
		\bibitem{Salia}
		N. Salia, P\'osa-type results for Berge hypergraphs, \emph{Electron. J. Combin.} 31 (2024), no.~2, Paper No.~P2.42.
		
		\bibitem{ClemensEhrenmullerPerson}
		D. Clemens, J. Ehrenm\"uller and Y. Person, A Dirac-type theorem for Berge cycles in random hypergraphs, \emph{Electron. J. Combin.} 27 (2020), no.~3, Paper No.~P3.39.
	
	\bibitem{ChenZhang}
A. Chen and L. Zhang, Asymptotically optimal local resilience of Berge
Hamilton cycles in random $k$-uniform hypergraphs, 
submitted.
	
	
	
		\bibitem{ChenZhang3}
		A. Chen and L. Zhang, Minimum degree condition of Berge Hamiltonicity
		in random $3$-uniform hypergraphs, \emph{Filomat} 37 (2023), no.~26,
		9039--9050.
		
	
\bibitem{BalBerkowitzDevlinSchacht}
D. Bal, R. Berkowitz, P. Devlin and M. Schacht, Hamiltonian Berge cycles in random hypergraphs, \emph{Combin. Probab. Comput.} 30 (2021), 228--238.
	
		
		\bibitem{AllenParczykPfenninger}
		P. Allen, O. Parczyk and V. Pfenninger, Resilience for tight Hamiltonicity, \emph{Combinatorial Theory} 4 (2024), no.~1, Paper No.~9.
		
		\bibitem{PetrovaTrujic}
		K. Petrova and M. Truji\'c, Transference for loose Hamilton cycles in random 3-uniform hypergraphs, \emph{Random Structures Algorithms} 65 (2024), 313--341.
	
		
		\bibitem{KeevashEtAl}
		P. Keevash, D. K\"uhn, R. Mycroft and D. Osthus, Loose Hamilton cycles in hypergraphs, \emph{Discrete Math.} 311 (2011), 544--559.
		
		\bibitem{NenadovSkoric}
		R. Nenadov and N. \v{S}kori\'c, Powers of Hamilton cycles in random graphs and tight Hamilton cycles in random hypergraphs, \emph{Random Structures Algorithms} 54 (2019), 187--208.
		
		\bibitem{AlonSpencer}
		N. Alon and J.~H. Spencer, \emph{The Probabilistic Method}, 4th ed., Wiley, Hoboken, NJ, 2016.
		
		\bibitem{Posa1976}
		L. P\'osa, Hamiltonian circuits in random graphs, \emph{Discrete Math.} 14 (1976), 359--364.
		
	\end{thebibliography}
\end{document}